\documentclass[article]{article}

\newcommand{\R}{\mathbb{R}}
\newcommand{\C}{\mathbb{C}}

\let \eps\varepsilon

\newcommand{\tmu}{\widetilde{\mu}}
\newcommand{\hmu}{\widehat{\mu}}
\newcommand{\hpsi}{\widehat{\psi}}
\newcommand{\tpsi}{\widetilde{\psi}}
\newtheorem{remark}{Remark}
\newtheorem{algorithm}{Algorithm}
\usepackage{epsfig} 
\usepackage{amsmath} 
\usepackage{amssymb}  

\title{\LARGE \bf
A greedy algorithm for the identification of quantum systems. 
}

\author{
Yvon Maday\thanks{Y. Maday, UPMC Univ Paris 06, UMR 7598, Laboratoire Jacques-Louis Lions, F-75005, Paris, France
 and Division of Applied Mathematics, Brown University, Providence, RI, USA. {\tt\small maday@ann.jussieu.fr}} 
and Julien Salomon\thanks{J. Salomon, CEREMADE, UMR CNRS 7534, Universit\'e Paris IX
Dauphine, Pl. de Lattre de Tassigny, 75775 Paris Cedex 16, France. {\tt\small  salomon@ceremade.dauphine.fr}}
}

\begin{document}

\maketitle
\thispagestyle{empty}
\pagestyle{empty}

\begin{abstract}
The control of quantum phenomena is a topic that has carried out many challenging
problems. Among others, the Hamiltonian identification, i.e, the
inverse problem associated with the unknown features of a quantum system
is still an open issue.
In this work, we present an algorithm that enables to design a set of
selective laser fields that can be used, in a second stage, to identify
unknown parameters of quantum systems.
\end{abstract}
\section{Introduction}\label{intro}
The possibility to use coherent light to manipulate molecular systems
at the nanoscale has been demonstrated both theoretically \cite{karine} and
experimentally \cite{weinacht}. Different types of methods have proven their
relevancy for various settings, ranging from electron to large
polyatomic molecules \cite{hb,hr0,rice,hr2,hr1}.\\
At the same time, the ability to generate a large amount of quantum
dynamics data in a small time frame can also be used to extract from
experiments the values of unknown parameters of quantum systems. The
corresponding inverse problem, usually called Hamiltonian
identification has recently been subject to significant developments through
encouraging experimental
results \cite{hr3}. \\
Various formulations in an optimization settings have been
studied. Because of the nature of the available data, zero order
methods were first tested, see e.g. the
technique of map inversion \cite{hr4}. The use of 
optimal control techniques was
then introduced \cite{mazyar,hr5}.\\
Contrary to this last class of methods, we present in this
work a methodology that enables to handle situations where the
experimental measurements are
provided only at a given time. Our approach is based on a precomputation
that provides a family of {\sl selective} laser fields. Roughly speaking,
these laser fields are
designed iteratively to highlight variations in the
parameters that are subject to the identification.  In
a second stage, these fields and the experimental measurements are
used to assemble a nonlinear system satisfied by the to-be-identified 
parameters. \\
The paper is organized as follows: the optimization framework and the
assumptions we use are presented in Sec. \ref{identifpb}. In
Sec. \ref{algo}, the structure of our algorithm is given. The
procedures used in the two parts of this algorithm are described in
Sec. \ref{equiv} and Sec. \ref{disc}. The identification step is
explained in Sec. \ref{idproc}. Details about practical implementation
and some numerical results are given Sec. \ref{test}. We conclude with
some remarks in Sec. \ref{conc}.\\
Throughout this paper, $\Omega$ is a spacial domain in $\R^d$, $d=1,2,3$,
$L^2$ denotes the space of complex valued square integrable functions over $\Omega$,  and $\langle ., .\rangle$
 the usual Hermitian product associated to $L^2$. The following
 standard convention is used:  
$$\langle a |O| b \rangle :=\langle a,O(b)\rangle, \ a\in L^2,\
b\in L^2,\ O\in{\cal L}(L^2;L^2)$$  the set of all linear operator from $L^2$ into $L^2$.
Finally, we use
 $\Re(z),\Im(z)$ to denote respectively the real and the imaginary
 part of a complex number $z$. \\

\section{The identification problem}\label{identifpb}
We first introduce the model and the framework
used in this paper. 
\subsection{Control of the Schr\"odinger Equation}
Consider a quantum system $\psi\in H^1$, with norm $\|\psi\|_{L^2} = 1$,  evolving according to the Schr\"odinger equation
\begin{equation}\label{schrod}
\left\{\begin{array}{ccl}
    i\dot \psi &=& [H_0+V+\eps(t)\mu]\psi\\
\psi(0)&=&\psi_0, 
\end{array}
\right.
\end{equation}
where $H_0$ is the kinetic energy operator, $V\in{\cal L}(L^2;L^2)$ the potential operator
and $\mu\in{\cal L}(L^2;L^2)$ the dipole moment operator coupling
the system to a time-dependent external laser field $\eps(t)$.  In this
context, $\eps$ reads as a control since it can be chosen by the experimenter.\\
 
In the settings we consider here, we assume that the internal Hamiltonian $H=H_0+V$ is
known so that the goal is to identify the dipole moment operator $\mu$. The generalization to the identification of $V$
should not give rise to any particular problem and is left to a future contribution.

The basic hypothesis made on $\mu$ is that it belongs to (or actually can be conveniently approximated by)  a finite dimensional space spanned by some basis set $\mathcal{B}_\mu=(\mu_\ell)_{\ell=1,\dots, L}$.

\subsection{Experimental measurements and controllability}
In order to perform the identification, we assume that
 given a time $T$ and a laser field $\eps\in L^2(0,T)$, the
 experimenter can measure, for some fixed state $\psi_1\in L^2$, with norm $\|\psi_1\|_{L^2} = 1$, the
 value $\varphi(\mu,\eps):=\langle \psi_1,\psi(T)\rangle$.\\
Note that all what follows still holds when considering several
measurements a time $T$, i.e., in the case where a set of measurement
$\left(\langle \psi_\ell,\psi(T)\rangle\right)_{\ell=1\dots, p}$, with
$p>1$, is known.\\
Finally, we assume that the system under consideration is wavefunction
controllable, i.e., that $\eps\in L^2(0,T)\mapsto \psi(T) $ is surjective.
\subsection{Formulation of problem}
Our identification method is based on a particular formulation of the
identification problem that we now briefly introduce.\\
Denote by $\mu^\star$ the actual dipole moment operator of a given system. The solution $\mu=\mu^\star$
of our problem also solves the minimization problem:
\begin{equation}\label{infsup}
\inf_{\mu \in {\cal L}(L^2;L^2)}\sup_{\eps \in L^2(0,T)} \vert
\varphi(\mu,\eps)-\varphi(\mu^\star,\eps) \vert^2.
\end{equation}
This settings highlights the fact that as long as $\mu\neq\mu^\star$, a
selective laser field should be designed so that the difference between
$\mu$ and $\mu^\star$ is discerned through the measurement
$\varphi(\mu,\eps)$. \\

\section{Structure of the algorithm}\label{algo}
Our algorithm consists in designing, through a finite iterative procedure, a set of selective laser
fields. We start with the general
structure of our algorithm. Details about its 
steps are given in the next sections.
\subsection{The selective laser fields computation greedy algorithm}
Starting from the  basis set
$\mathcal{B}_\mu=(\mu_\ell)_{\ell=1,\dots,L}$, the algorithm builds up iteratively a
set of $L$ selective laser fields as follows. \\
\begin{algorithm}\label{ga}(Selective laser fields computation greedy algorithm)
Let us define $\eps^1$ a laser field that solves the
problem:
$$\sup_{\eps \in L^2(0,T)} \vert
\varphi(\mu_1,\eps)
\vert^2.$$
Suppose now that at the step $k$, with $1<k\leq L$, a laser field $\eps^{k-1}$ is
given. The computation of $\eps^k$ is
performed according to the two following sub-steps:

\begin{enumerate}
\item \underline{Fitting step} : Find $(\alpha^{k}_j)_{j=1,\dots,k-1}$ that solves the problem:
\begin{equation}\label{minsq1}
\left\{\begin{array}{ccl}
\varphi(\sum_{j=1}^{k-1}\alpha^k_j\mu^j,\eps^1)&=&\varphi(\mu^k,\eps^1)\\
&\vdots & \\
\varphi(\sum_{j=1}^{k-1}\alpha^k_j\mu^j,\eps^m)&=&\varphi(\mu^k,\eps^m)\\
&\vdots & \\
\varphi(\sum_{j=1}^{k-1}\alpha^k_j\mu^j,\eps^{k-1})&=&\varphi(\mu^k,\eps^{k-1}),
\end{array}\right.
\end{equation}
in the minimum mean square error sense.
\item \underline{Discriminatory step} : Find $\eps^k$ that solves the problem:
$$\eps^k = \hbox{argmax}_{\eps \in L^2(0,T)} \vert \varphi(\mu^k,\eps)-\varphi(\sum_{j=1}^{k-1}\alpha^k_j\mu^j,\eps)
\vert^2.$$
\end{enumerate}
\end{algorithm}

The initialization of the algorithm is somehow arbitrary, the only
requirement is that $\eps^1$ has a link with the type of
measurement. In our case, we decide to maximize it.
\begin{remark}
Note that, in opposition to usual approaches (see e.g. \cite{hr5,mazyar}), our method
plays the role of a precomputation step since the actual measurements
$\varphi(\mu^\star,\eps)$ are not required at this stage. 
\end{remark}
\subsection{Intuitive interpretation of the algorithm}
In the first sub-step of an iteration of Algorithm \ref{ga}, one looks for a defect
of selectivity of the
current laser fields $\eps^1,\dots,\eps^{k-1}$: in the case the minimum reaches zero, two
distinct dipole moment operators give rise to two identical measurements when exited with the laser fields $\eps^1,\dots,\eps^{k-1}$.
On the contrary, the second sub-step aims at computing a laser field that
compensates this defect. 
 These two sub-steps
corresponds respectively to the minimization part and to the
maximization part of the formulation \eqref{infsup}.
\begin{remark}
Even if no hierarchy is assumed in the basis $\mathcal{B}_\mu$, this
algorithm should be viewed as a first step towards future works that
handle infinite dimensional systems. In such a framework, the sum
$\sum_{j=1}^{k-1}\alpha^k_j\mu^j$ would read as an asymptotic
expansion of the dipole moment operator. 
\end{remark}

This algorithm belongs to the class of greedy algorithms, since it follows the problem-solving's heuristic of making the locally optimal choice (in the second sub-step) at each stage
with the hope of finding the global optimum that solves \eqref{infsup}.\\

\section{Fitting step}\label{equiv}
Let us first focus on the first sub-step of the algorithm. Consider an
integer $k$ such that $1<k\leq L$ and denote by
$K^k$ the functional (defined on $\R^{k-1}$):
$$K^k(\alpha)=\sum_{m=1}^{k-1}\vert
\varphi(\mu^k,\eps^m)-\varphi(\sum_{j=1}^{k-1}\alpha_j\mu^j,\eps^m)
\vert^2.$$
During this sub-step, one has to find the minimum of the cost
functional $K^k$. To do this, a standard global minimization
algorithm applied to this minimum mean square error associated problem.\\
Note that, for small values of $L$, the gradient of the functional $K^k$ can
be computed thanks to the formula:\\\\
$\nabla K^k(\alpha).\delta \alpha=$ 
$$ \sum_{m=1}^{k-1}2\Re\left( \langle\psi_{\eps^m}^\alpha(T)-\psi_{\eps^m}^k(T),\psi_1\rangle\langle \delta \psi_{\eps^m}^\alpha (T),\psi_1\rangle\right), $$ 
where $\psi_{\eps^m}^\alpha$ and  $\psi_{\eps^m}^k$ are the solutions of Eq. \eqref{schrod}
with $\eps=\eps^m$ as laser field, and
$\mu=\sum_{j=1}^{k-1}\alpha_j\mu^j$ and $\mu=\mu^k$ respectively as
dipole moment operator. The variation $\delta \psi^\alpha$ is computed
thanks to:
\begin{equation}\nonumber
\left\{\begin{array}{ccl}
    i \delta\dot \psi_{\eps^m}^\alpha &=& \eps^{k-1}\left(\sum_{j=1}^{k-1}\alpha_j\mu^j\right)
    \delta \psi_{\eps^m}^\alpha \\&&+
    [H_0+V+\eps^{k-1}(t)\left(\sum_{j=1}^{k-1}\delta
    \alpha_j\mu^j\right)]\psi_{\eps^m}^\alpha\\
\delta \psi_{\eps^m}^\alpha(0)&=&0. 
\end{array}
\right.
\end{equation}
In this way the computation of the components of $\nabla K^k(\alpha)$ can be parallelized to make the use of
gradient methods feasible.\\

\section{Discriminatory step}\label{disc}
To achieve the second sub-step of Algorithm \ref{ga}, we adapt an efficient strategy
usually used in in quantum control. This strategy has given rise to a
large class of algorithms often called "monotonic schemes". For a
general presentation of these algorithms, we refer to \cite{jsgt}.
\subsection{Improvement of the selectivity of a given laser field}\label{discms}
Let us present in more details how this strategy applies in our case.
Note first that, given a laser field $\eps\in L^2(0,T)$, and two dipole
moment operators $\tmu$ and $\hmu$, one has:
$$\vert
\varphi(\tmu,\eps)-\varphi(\hmu,\eps)
\vert^2=\langle \tpsi(T)-\hpsi(T) |O_{\psi_1}|
\tpsi(T)-\hpsi(T)  \rangle,$$
 where $O_{\psi_1} = \psi_1.\psi_1^T$, $\tpsi$ and $\hpsi$ are the solutions of
Eq. \eqref{schrod} with respectively $\mu=\tmu$ and $\mu=\hmu$ as dipole
moment operator.  

In order to compare the selectivity of $\eps$ and $\eps'$, we introduce the
functional:
$$J(\eps)=\langle \tpsi(T)-\hpsi(T) |O_{\psi_1}|
\tpsi(T)-\hpsi(T)  \rangle-\beta\int_0^T\eps^2(t)dt,$$
which has to be maximized. For sake of simplicity, we omit the
dependence of $J$ with $\tmu$ and $\hmu$ in the notations.\\
The additional term
$\beta\int_0^T\eps^2(t)dt$, is introduced for two complementary reasons: first, as
it penalizes the $L^2$-norm of the laser field, it
enables to obtain feasible laser fields and secondly, it improves the
convergence of Algorithm \ref{ms} below.\\
Consider now another laser field $\eps'\in L^2(0,T)$, and denote by
$\tpsi'$ and $\hpsi'$ the corresponding solutions of
Eq. \eqref{schrod} with $\mu=\tmu$ and $\mu=\hmu$ respectively. We
introduce the two adjoints states defined by:
\begin{equation}
\left\{
\begin{array}{ccl}
i\dot{\widetilde{\chi}}   &=& [H_0+V+\eps(t)\tmu]\widetilde{\chi}\\
      \widetilde{\chi}(T) &=&O_{\psi_1}\left(\tpsi(T)-\hpsi(T)\right),
\end{array}
\right.\label{adj1}
\end{equation}
and
\begin{equation}
\left\{\begin{array}{ccl}
i\dot{\widehat{\chi}    } &=& [H_0+V+\eps(t)\hmu]\widehat{\chi}\\
      \widehat{\chi}(T)   &=& O_{\psi_1}\left(\tpsi(T)-\hpsi(T)\right).
\end{array}
\right.\label{adj2}
\end{equation}
One has:
\begin{eqnarray}
J(\eps')-J(\eps)=\langle \delta \psi'(T)-\delta \psi(T) |O_{\psi_1}|
\delta \psi'(T)-\delta \psi(T) \rangle\nonumber\\
+2\Re \langle \delta \psi'(T)-\delta \psi(T),\widetilde{\chi}(T)-\widehat{\chi}(T) \rangle
\nonumber\\
-\beta\int_0^T \eps'^2(t)-\eps^2(t)dt
\nonumber\\
=\langle \delta \psi'(T)-\delta \psi(T) |O_{\psi_1}|
\delta \psi'(T)-\delta \psi(T) \rangle\nonumber\\
  + \int_0^T\left(\eps'(t)-\eps(t)\right)\nonumber\\
          \left( 2\Im\langle
            \widetilde{\chi}(t)|\tmu|\tpsi'(t)\rangle\!  - \!
                 2\Im\langle \widehat  {\chi}(t)|\hmu|\hpsi'(t)\rangle  
\! -\!\beta\left(\eps'(t)\! +\! \eps(t)\right)\right)dt,\nonumber
\end{eqnarray}
\begin{equation}
\label{identity}
\end{equation}
where we denote $\delta \psi'(T)=\tpsi'(T)-\hpsi'(T)$ and $\delta \psi(T)=\tpsi(T)-\hpsi(T)$.
Identity \eqref{identity} gives a criterion to guarantee that
$\eps'$ is more selective than $\eps$. Indeed, suppose that $\eps'$ satisfies for all
$t\in [0,T]$ the condition:
\begin{eqnarray}
\left(\eps'(t)-\eps(t)\right)   
\left( 2\Im\langle \widetilde{\chi}(t)|\tmu|\tpsi'(t)\rangle  -
       2\Im\langle \widehat  {\chi}(t)|\hmu|\hpsi'(t)\rangle \right.\nonumber\\
\left. \phantom{\Im\langle \widetilde{\chi}(t)|\tmu|\tpsi'(t)\rangle} -\beta(\eps'(t)+\eps(t))\right)\geq 0,\label{criterion}
\end{eqnarray}
then $J(\eps')\geq J(\eps)$.\\
Various ways to ensure that \eqref{criterion} holds. For example \cite{tannor}, one can define
$\eps'$ at each time $t$ as the solution of the equation:
\begin{eqnarray}
\eps'(t)-\eps(t)= \frac\theta\beta\left(
2\Im\langle \widetilde{\chi}(t)|\tmu|\tpsi'(t)\rangle  -
2\Im\langle \widehat  {\chi}(t)|\hmu|\hpsi'(t)\rangle \right. \nonumber\\
 \left. \phantom{\Im\langle \widetilde{\chi}(t)|\tmu|\tpsi'(t)\rangle}
   -\beta\left(\eps'(t)+\eps(t)\right)\right),\nonumber
\end{eqnarray}
\begin{equation}
\label{thetastrat}
\end{equation}
where $\theta$ is a given strictly positive number. 
In this case, one has:
\begin{eqnarray}
J(\eps')-J(\eps)&\!\!=\!\!&\langle \delta\psi'(T)-\delta\psi(T)
|O_{\psi_1}|
\delta\psi'(T) - \delta\psi(T)\rangle\nonumber\\
&& +\frac\beta\theta \int_0^T\left(\eps'(t)-\eps(t)\right)^2dt\geq 0,\nonumber
\end{eqnarray}
which is the desired conclusion. In Sec. \ref{ns}, we present an
alternative that can be obtained in a time discretized settings.
\subsection{Discriminatory sub-algorithm}
We derive form the previous considerations the following iterative procedure
to define a laser field $\eps^{k}$ that maximizes $J(\eps)$:\\
\begin{algorithm}\label{ms}(Discriminatory sub-algorithm)
Let  $Tol$ be a positive number. Consider an initial guess $\eps^k_0$ and compute the corresponding
solutions of Eq. \eqref{schrod} with $\tmu$ and $\hmu$, say $\tpsi_0$ and
$\hpsi_0$. Set $err=2.Tol$.\\
While $err>Tol$, do:
\begin{enumerate}
\item  Use Eqs. (\ref{adj1}--\ref{adj2}) with $\eps=\eps^k_\ell$, $\tpsi=\tpsi_\ell$ and
  $\hpsi=\hpsi_\ell$, to compute $\widetilde{\chi}_\ell $ and $\widehat{\chi}_\ell$, respectively.
\item Compute simultaneously the laser field $\eps^k_{\ell+1}$ and the
  states $\tpsi_{\ell+1}$ and $\hpsi_{\ell+1}$ the solutions of
  coupled system composed of Eq.
  \eqref{thetastrat} with $\widetilde{\chi}= \widetilde{\chi}_\ell$,
  $\widehat{\chi}= \widehat{\chi}_\ell$ and Eq.
  \eqref{schrod} with $\mu=\tmu$ and $\mu=\hmu$ respectively.
\item $\ell\leftarrow \ell+1$, $err=|\eps^k_{\ell+1}-\eps^k_{\ell}|$.
\end{enumerate}
\end{algorithm}

In \cite{mst}, one shows that Eq. \eqref{thetastrat} has a solution and presents some efficient numerical nonlinear
solvers to compute it.\\

\section{Identification procedure}\label{idproc}

Once the $L$ selective fields $\eps^1,...,\eps^L$ have been computed,
one can use them experimentally to obtain the corresponding
measurements
$\varphi(\mu^\star,\eps^1),...,\varphi(\mu^\star,\eps^L)$. \\
The identification procedure consists then in finding the linear
combination $(\alpha^1,...,\alpha^L)$ that solves the following
nonlinear system:
\begin{equation}\label{minsq2}
\left\{\begin{array}{ccl}
\varphi(\sum_{j=1}^{L}\alpha_j\mu^j,\eps^1)&=&\varphi(\mu^\star,\eps^1)\\
&\vdots & \\
\varphi(\sum_{j=1}^{L}\alpha_j\mu^j,\eps^k)&=&\varphi(\mu^\star,\eps^k)\\
&\vdots & \\
\varphi(\sum_{j=1}^{L}\alpha_j\mu^j,\eps^L)&=&\varphi(\mu^\star,\eps^L).
\end{array}\right.
\end{equation}
in the mean square sense.
In this view, the standard global optimization procedure used for the first sub-step
of algorithm can be applied to the 
associated problem.\\
Note that, in a finite-dimensional settings, the existence of
a solution is guaranteed.\\

\section{Numerical implementation and results}\label{test}
We give here details about the practical implementation of
 Algorithm \ref{ga}, 
and show its efficiency on an example.
\subsection{Numerical solvers}\label{ns}
In order to solve numerically Eq. \eqref{schrod}, we use the second
order Strang operator splitting \cite{strang}. Given $M>0$, a time step $\Delta t$
such that $M.\Delta t=T$ and an approximation $\psi_j$ of
$\psi(j.\Delta t)$ with $j<M$, this method leads in our case to
the following iteration:
\begin{equation}\label{strang}
\psi_{j+1}=e^{i H \frac{\Delta t}2}e^{i \eps_j \mu\Delta t}e^{i H
  \frac{\Delta t}2}\psi_j.
\end{equation}
In the second sub-step of Algorithm \ref{ga}, 
Discriminatory sub-algorithm \ref{ms} is adapted to this discrete settings. In this way, we
consider the time-discretized version of the cost functional $J$:
$$J_{\Delta t}(\eps)=\langle \tpsi_M-\hpsi_M |O_{\psi_1}|
\tpsi_M-\hpsi_M  \rangle-\beta\Delta t\sum_{j=0}^{M-1}\eps_j^2,$$
where  $\eps\in\R^{M-1}$. Fix now two discrete laser fields $\eps$ and $\eps'$, one can then repeat the computation done in
Sec. \ref{discms} to obtain:
\begin{eqnarray}
J_{\Delta t}(\eps')-J_{\Delta t}(\eps)=\langle \delta \psi'_M-\delta \psi_M |O_{\psi_1}|
\delta \psi'_M-\delta \psi_M  \rangle \nonumber\\
  +  \Delta t \sum_{j=0}^{M-1} \left( \eps'_j-\eps_j \right)\nonumber\\
          \left( 2\Im\langle \widetilde{\chi}_j|\tmu_{\Delta t}(\eps'_j,\eps_j)|\tpsi'_j\rangle
            -
                 2\Im\langle \widehat  {\chi}_j|\hmu_{\Delta
                   t}(\eps'_j,\eps_j)|\hpsi'_j\rangle \right. \nonumber\\
\left.
\phantom{ 2\Im\langle \widetilde{\chi}_j |\hmu_{\Delta
                   t}(\eps'_j,\eps_j)|  \hpsi'_j\rangle}  
-\beta\left(\eps'_j+\eps_j\right)\right),\label{identity2}
\end{eqnarray}
where the vectors $\widetilde{\chi}$, $\widehat{\chi}$, $\tpsi'$ and
$\hpsi'$ are computed using the iteration \eqref{strang} with $\mu=\tmu_{\Delta t}(\eps'_j,\eps_j)$ and $\mu=\hmu_{\Delta
  t}(\eps'_j,\eps_j)$. These matrices are the approximations of $\tmu$ and $\hmu$ respectively defined by:
$$
\begin{array}{ccl}
\tmu_{\Delta t}(\eps'_j,\eps_j) &=&e^{-i H\frac{\Delta
    t}{2}}\dfrac{e^{i\eps'_j \tmu \Delta t}-e^{i\eps_j \tmu \Delta
    t}}{i\Delta t(\eps'_j-\eps_j)}e^{i H\frac{\Delta t}{2}}\\
\hmu_{\Delta t}(\eps'_j,\eps_j) &=&e^{-i H\frac{\Delta
    t}{2}}\dfrac{e^{i\eps'_j \hmu \Delta t}-e^{i\eps_j \hmu \Delta
    t}}{i\Delta t(\eps'_j-\eps_j)}e^{i H\frac{\Delta t}{2}}.
\end{array}
$$

For the sake of simplicity, instead of solving the discrete version of
Eq. \eqref{thetastrat}, we compute $\eps'_j$ using one step of a Newton
optimization method applied to its corresponding term in the sum of
Eq. \eqref{identity2}. This strategy, and the one corresponding to Eq.
\eqref{thetastrat} are presented in more details in
\cite{mst}. Their convergence are proven in \cite{salomon}.

\subsection{Numerical test}
\subsubsection{Settings}
To illustrate the ability of our approach, we consider a simple finite dimensional settings where $H_0,\ V$ and $\mu$ are $3\times 3$ Hermitian matrices with
entries in $\C$ and $\psi(t)\in\C^3$. The internal Hamiltonian we consider is:
$$H=10^{-2}\left(
\begin{array}{ccc}
    1 &   0 &   0\\
    0 &   2 &   0\\
    0 &   0 &   4
\end{array}\right).
$$
Since Eq. \eqref{schrod} with such an internal Hamiltonian is generically controllable, we choose to
define the basis $\mathcal{B}_\mu$ randomly so that the systems
handled by our algorithm are almost surely controllable.\\
In order to work in a general framework, we chose $\mu^\star$
also randomly. In our example, we consider:

$$\mu^\star =\left(
\begin{array}{ccc}
  2.4154 &   1.9335  &  1.5822\\
  1.9335 &   1.4366  &  1.5991\\
  1.5822 &   1.5991  &  1.9843
\end{array}\right).
$$
The states $\psi_0$ and $\psi_1$ are 
$$\psi_0=\left(
\begin{array}{c}
1\\0\\0
\end{array}\right),\quad \psi_1=\left(
\begin{array}{c}
0\\0\\1
\end{array}\right).
$$
We choose $T=4000\pi$, which corresponds to $20$ periods of the
transition associated to the smallest frequency of 
the system.
\subsubsection{Algorithm parameters}
The minimum mean square error problems \eqref{minsq1}-\eqref{minsq2} are solved by standard
pseudo-Newton solvers. In order to make a global search, we repeat the
minimization 10 times with random initialization. The parameter
$\beta$ is adapted to make Algorithm \ref{ms} converge. In our
case $\beta=10^{-2}$.
\subsubsection{Numerical results}
\begin{figure}[htbp]
  \begin{center}
\includegraphics[height=5.5cm,width=\linewidth]{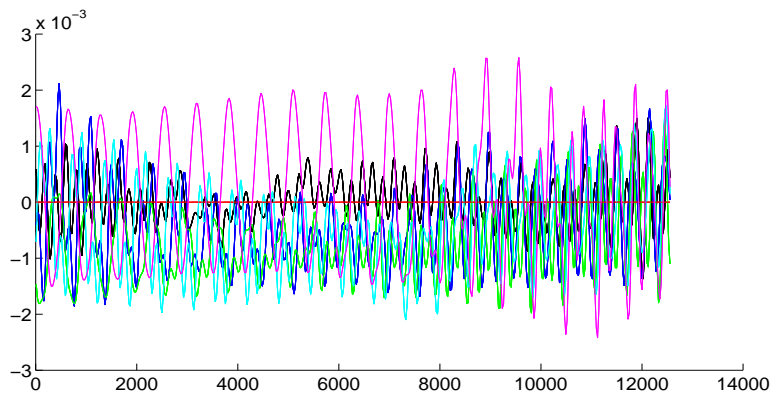}
    \caption{Selective laser fields obtained by Algorithm \ref{ms}.} \label{fig1}
  \end{center}
\end{figure}

The precomputation is achieved by our algorithm in approximately 80
min CPU.
The dipole moment operator is regained with a relative error 
$$\frac{\|\mu^\star-\mu\|_2}{\|\mu^\star\|_2}\approx   9.8960e-04,$$
in approximately 10 min CPU. The selective fields that have been obtained are depicted in Fig. \ref{fig1}.

\section{Concluding remarks}\label{conc}
The Selective laser fields computation greedy algorithm 
presented in this paper shows a good efficiency in a general
settings. However, there is some room for improvement of our
strategy. First, the choice of the basis $\mathcal{B}_\mu$ could
be improved, e.g. through an iterative procedure. Secondly, the
experimental measurements could be used during the computation of the
selective fields in order to design an online procedure. Lastly, some
work has to be done to design a more 
specific approach to treat the first sub-step of the algorithm. The
identification procedure presented in Sec. \ref{idproc} would also
certainly take advantage of such a study.

\section{ACKNOWLEDGMENTS}
The problem of identification in this context was raised during discussions with H. Rabitz from Princeton University and G. Turinici from Dauphine University, we thank them for helpful inputs.
This work was supported by the french A.N.R, ``Programme blanc C-Quid'' and PICS CNRS-NSF collaboration between the Department of Chemistry, Princeton University, and University Paris Dauphine.


\end{document}